\documentclass[a4paper]{amsart}

\usepackage{amssymb,amsmath,graphics}
\usepackage{latexsym}
\usepackage{eucal}
\usepackage[dvips]{graphicx}
\usepackage{color}
\usepackage{enumerate}
\makeatletter
\@addtoreset{equation}{section}\makeatother

\newtheorem{theo}{Theorem}[section]
\newtheorem{lem}[theo]{Lemma}
\newtheorem{prop}[theo]{Proposition}
\newtheorem{cor}[theo]{Corollary}

\newcommand{\mc}{\mathcal}
\newcommand{\rr}{\mathbb{R}}
\newcommand{\nn}{\mathbb{N}}
\newcommand{\cc}{\mathbb{C}}
\newcommand{\hh}{\mathbb{H}}
\newcommand{\zz}{\mathbb{Z}}

\newcommand{\la}{\lambda}
\newcommand{\eps}{\epsilon}

\newcommand{\pl}{\partial}
\newcommand{\x}{\times}

\newcommand{\cjd}{\rangle}
\newcommand{\cjg}{\langle}

\newcommand{\demi}{\frac{1}{2}}
\newcommand{\ndemi}{\frac{n}{2}}
\newcommand{\tra}{\textrm{Tr}}

\newcommand{\rang}{\textrm{rank}}

\newcommand{\indic}{\operatorname{1\negthinspace l}}

\renewcommand{\Re}{{\rm Re}}
\renewcommand{\Im}{{\rm Im}}
\newcommand{\R}{\mathbb{R}}
\newcommand{\C}{\mathbb{C}}

\renewcommand{\H}{\mathbb{H}^2}
\def\qed{\hfill$\square$}
\begin{document}
\title{Wave decay on convex co-compact hyperbolic manifolds}

\author[Colin Guillarmou]{Colin Guillarmou}
\address{Laboratoire de Mathematiques J. Dieudonn\'e \\
         Universit\'e de Nice, Parc Valrose\\
        Nice, France}
     \email{cguillar@math.unice.fr}
     
\author[Fr\'ed\'eric Naud]{Fr\'ed\'eric Naud}
\address{Laboratoire d'Analyse non-lin\'eaire et g\'eom\'etrie\\
Universit\'e d'Avignon\\
33 rue Louis Pasteur, 84000 Avignon\\
France}
\email{naud@univ-avignon.fr}

%

\begin{abstract}
\noindent  
For convex co-compact hyperbolic quotients $X=\Gamma\backslash\hh^{n+1}$, we analyze the long-time asymptotic 
of the solution of the wave equation $u(t)$ with smooth compactly supported initial data $f=(f_0,f_1)$. 
We show that, if the Hausdorff dimension $\delta$ of the limit set is less than $n/2$, then 
$u(t)=C_\delta(f)e^{(\delta-\ndemi)t}/\Gamma(\delta-n/2+1)+e^{(\delta-\ndemi)t}R(t)$ where 
$C_{\delta}(f)\in C^\infty(X)$ and $||R(t)||=\mc{O}(t^{-\infty})$. 
We explain, in terms of conformal theory of the conformal infinity of $X$, the special cases $\delta\in n/2-\nn$ where the leading asymptotic term vanishes. In a second part, we show for all $\eps>0$ the existence of an infinite number of resonances (and thus zeros of Selberg zeta function) in the strip $\{-n\delta-\eps<\Re(\la)<\delta\}$. As a byproduct 
we obtain a lower bound on the remainder $R(t)$ for generic initial data $f$.
\end{abstract}

\maketitle

\section{Introduction}
It is well-known that on a compact Riemannian manifold $(X,g)$, any solution $u(t,z)$ of the wave equation $(\pl_t^2+\Delta_g)u(t,z)=0$ expands 
as a sum of oscillating terms of the form $e^{i\la_jt}a_j(z)$ where $\la_j^2$ are 
the eigenvalues of the Laplacian $\Delta_g$ and $a_j$ some associated eigenvectors. The eigenvalues then give the 
frequencies of oscillation in time. For non-compact manifolds, the situation is much more complicated and 
no general theory describes the behaviour of waves as time goes to infinity, at least in terms of 
spectral data. A first satisfactory description has been given by Lax-Phillips \cite{LP} and Vainberg \cite{V} 
for the Laplacian $\Delta_{X}$ with Dirichlet condition on $X:=\rr^n\setminus\mc{O}$ where $\mc{O}$ is a compact 
obstacle and $n$ odd; indeed if $u(t)$ is the solution of $(-\pl_t^2-\Delta_X)u(t,z)=0$ with compactly supported smooth initial data in $X$ and under a \emph{non-trapping} condition, they show an expansion as $t\to+\infty$ of the form 
\[u(t,z)=\sum_{\substack{\la_j\in\mc{R}\\
\Im(\la_j)<N}}\sum_{k=1}^{m(\la_j)}e^{i\la_jt}t^{k-1}u_{j,k}(z)+\mc{O}(e^{-(N-\eps)t}), \quad \forall N>0, \forall \eps>0\]
uniformly on compacts, where $\mc{R}\subset \{\la\in\cc,\Im(\la)\geq 0\}$ 
is a discrete set of complex numbers called \emph{resonances} 
associated with a multiplicity function $m:\mc{R}\to \nn$, and $u_{j,k}$ are smooth functions.  
The real part of $\la_j$ is a frequency of oscillation while the imaginary part is an exponential decay rate of the 
solution. Resonances can in general be defined as poles of the meromorphic continuation of the 
Schwartz kernel of the resolvent of $\Delta_X$ through the continuous spectrum. 
   
In \cite{TZ}, Tang and Zworski extended this result for \emph{non-trapping} black-box perturbation of $\rr^n$ and 
considered also a strongly trapped setting, namely when
there exist resonances $\la_j$ such that\footnote{This 
is typically the case when $P$ has elliptic trapped orbits as shown in \cite{P}} $\Im(\la_j)<(1+|\la_j|)^{-N}$ for all $N>0$, satisfying in addition 
some separation and multiplicity conditions. 
The expansion of wave solutions then involved these resonances and the error is $\mc{O}(t^{-N})$ for all $N>0$.  
This last result has also been generalized by Burq-Zworski \cite{BZ} for semi-classical problems.

It is important to notice that such results are almost certainly not optimal when 
the trapping is hyperbolic since, at least for all known examples, resonances do not seem 
to approach the real line faster than polynomially. 
Christiansen and Zworski \cite{CZ} studied two examples in hyperbolic geometry, the modular surface and 
the infinite volume cylinder, they showed a full expansion of waves in terms of resonances 
with exponentially decaying error terms. The proof is based on a separation of variables computation in the cylinder case (here the trapping geometry is that of a single closed hyperbolic orbit)
while it relies on well-known number theoretic estimates for the Eisenstein series in the modular case. 
The case of De Sitter-Schwarzchild metrics has recently been studied by Bony-H\"afner \cite{BH}
using also separation of variables and rotational symmetry of the space. This is another example of hyperbolic trapping.
Clearly, the general hyperbolic trapping situation is an issue and the above results are always based on very explicit computations or the arithmetic nature of the manifold.
It is therefore of interest to consider more general cases of hyperbolic trapping geometries, the most basic examples being the convex co-compact quotients of the hyperbolic space $\hh^{n+1}$ that can be considered as the simplest non-trivial models of open quantum chaotic systems.
\\ 

Hyperbolic quotients $\Gamma\backslash\hh^{n+1}$ by a discrete group of isometries with
only hyperbolic elements (those that do not fix points in $\hh^{n+1}$ but fix two points on the sphere 
at infinity $S^n=\pl\hh^{n+1}$) and admitting a finite sided fundamental domain are called
\emph{convex co-compact}. The Laplacian on such a quotient $X$ has for continuous and essential 
spectrum the half-line $[n^2/4,\infty)$, the natural wave equation is 
\begin{equation}\label{waveeq}
(\pl_t^2+\Delta_X-n^2/4)u(t,z)=0, \quad u(0,z)=f_0(z), \quad \pl_tu(0,z)=f_1(z),
\end{equation}    
its solution is 
\begin{equation}\label{utf}
u(t)=\cos\Big(t\sqrt{\Delta_X-\frac{n^2}{4}}\Big)f_0+
\frac{\sin \Big(t\sqrt{\Delta_X-\frac{n^2}{4}}\Big)}{\sqrt{\Delta_X-\frac{n^2}{4}}}f_1.\end{equation}
For a convex co-compact quotient $X=\Gamma\backslash \hh^{n+1}$, 
the group $\Gamma$ acts on $\hh^{n+1}$ as isometries but also on the sphere at infinity 
$S^n=\pl\hh^{n+1}$ as conformal transformations. The limit set $\Lambda(\Gamma)$ 
of the group is  the set of accumulation points on $S^n$ of the orbit $\Gamma.m$ for the Euclidean
topology of the ball $\{z\in\rr^{n+1};|z|\leq 1\}$ for any picked $m\in\hh^{n+1}$, it is well known that $\Lambda(\Gamma)$ 
does not depend on the choice of $m$.
We denote by $\delta\in(0,n)$ the Hausdorff dimension of $\Lambda(\Gamma)$,
\[\delta:=\dim_{H}(\Lambda(\Gamma)).\] 
It is proved by Patterson \cite{Pat} and Sullivan \cite{SU} that $\delta$ is also the exponent of convergence of Poincar\'e series
\begin{equation}\label{poincareseries}
P_{\la}(m,m'):=\sum_{\gamma\in\Gamma}e^{-\la d_{h}(m,\gamma m')}, \quad m,m'\in\hh^{n+1},
\end{equation}   
where $d_h$ is the hyperbolic distance. Standard coordinates on the unit sphere bundle $SX=\{(z,\xi)\in TX; |\xi|=1\}$ show that $2\delta+1$ is the Hausdorff dimension of the trapped set of the geodesic flow on $SX$.

We denote by $\Omega:=S^n\setminus \Lambda(\Gamma)$ the domain of discontinuity of $\Gamma$,
this is the largest open subset of $S^n$ on which $\Gamma$ acts properly discontinuously.
The quotient $\Gamma\backslash\Omega$ is a compact manifold and $X$ can be compactified
into a smooth manifold with boundary $\bar{X}=X\cup \pl\bar{X}$ with $\pl\bar{X}=
\Gamma\backslash\Omega$. It turns out that $\pl\bar{X}$ inherits from the hyperbolic metric $g$
on $X$ a conformal class of metrics $[h_0]$, namely the conformal class of $h_0=x^2g|_{T\pl\bar{X}}$
where $x$ is any smooth boundary defining function of $\pl\bar{X}$ in $\bar{X}$.\\

In this paper, we focus on the case when $\delta<n/2$ since if $\delta>n/2$, the Laplacian $\Delta_X$ has pure point spectrum
in $(0,n^2/4)$ that gives the leading asymptotic behaviour of $u(t)$ by usual spectral theory. 
We prove the following result.
\begin{theo}\label{mainth0}
Let $X$ be an $(n+1)$-dimensional convex co-compact hyperbolic manifold such that $\delta<n/2$,
and let $f_0,f_1,\chi\in C_0^\infty(X)$. With $u(t)$ defined by \eqref{utf}, 
as $t\to +\infty$, we have the asymptotic
\begin{equation}\label{chiut}
\chi u(t)=\frac{A_X}{\Gamma(\delta-\ndemi+1)}e^{-t(\ndemi-\delta)}\cjg u_\delta,(\delta-\ndemi)f_0+f_1\cjd \chi u_\delta +\mc{O}_{L^2}(e^{(\delta-\ndemi)t}t^{-\infty})\end{equation}
where $u_\delta$ is the Patterson generalized eigenfunction defined in \eqref{udelta}
and $\cjg,\cjd$ is the distributional pairing, 
$A_X\in\cc\setminus\{0\}$ is a constant depending on $X$.
 \end{theo}

\textsl{Remark 1}: when $\delta\notin n/2-\nn$, this shows that the ``dynamical dimension'' $\delta$ controls the exponential decay rate of waves, or \emph{quantum decay rate}\footnote{This kind of result was predicted
in \cite{N2}.}. It seems to be the first rather general example of hyperbolic trapping for which we have an explicit 
asymptotic for the waves, in terms of geometric data. However, we point out that the recent work of Petkov-Stoyanov \cite{PeS} should in principle imply an expansion in terms of a finite number of resonances for the exterior problem with strictly convex obstacles. We also believe that a result similar to Theorem \ref{mainth} holds for general negatively curved 
asymptotically hyperbolic manifolds, this will be studied in a subsequent work.\\

\textsl{Remark 2}: 
In the special case $\delta\in n/2-\nn$ (note that it can happen only for $n\geq 3$ i.e. for four and higher dimensional manifolds)  
the leading term vanishes in view of the Euler $\Gamma$ fonction in \eqref{chiut}. Waves for 
this special case turn out to decrease faster. We explain this fact in the last section 
of the paper, and it is somehow related to the conformal theory of $\pl\bar{X}$: what happens is 
that when $\delta\notin n/2-\nn$, $\la=\delta$ is always the closest pole to the continuous spectrum of the 
meromorphic extension of the resolvent $R(\la):=(\Delta_X-\la(n-\la))^{-1}$ 
and $u_\delta$ is an associated non-$L^2$ eigenstate $(\Delta_X-\delta(n-\delta))u_\delta=0$, 
while when $\delta\in n/2-k$ with $k\in\nn$, the extended resolvent $R(\la)$ is holomorphic at $\la=\delta$ and
$u_\delta$ has asymptotic behaviour near $\pl\bar{X}$
\[u_\delta(z)= x(z)^{\delta}f_\delta +\mc{O}(x(z)^{\delta+1})\]
where $f_\delta\in C^\infty(\pl\bar{X})$ is an element of $\ker(P_k)$, $P_k$ being the $k$-th GJMS conformal Laplacian
\cite{GJMS} of the conformal boundary  $(\pl\bar{X},[h_0])$; more precisely $P_j>0$ for all $j=1,\dots,k-1$ while 
$\ker P_k=\textrm{Span}(f_\delta)$. The manifold has a special conformal geometry at infinity that makes the resonance $\delta$ disappears and transforms into a $0$-eigenvalue for the conformal Laplacian $P_k$.\\

The proof uses methods of Tang-Zworski \cite{TZ} together with informations on the 
closest resonance to the critical line, that is $\delta$ when $\delta\notin n/2-\nn$ 
(the physical sheet for the resolvent $R(\la):=(\Delta-\la(n-\la))^{-1}$
is $\{\Re(\la)>n/2\}$)
this last fact has been proved by Patterson \cite{Pa} using Poincar\'e series and Patterson-Sullivan 
measure. The powerful dynamical theory of Dolgopyat \cite{Do} has
been used by the second author \cite{N} (for surfaces) 
and Stoyanov \cite{S} (in higher dimension) to prove the existence of a strip with no zero
on the left of the first zero $\la=\delta$ for the Selberg zeta function. Using results
of Patterson-Perry \cite{PP}, this implies a strip $\{\delta-\eps<\Re(\la)<\delta)\}$ with no resonance. Then we can view $u(t)$ as a contour integral of the resolvent 
$R(\la)$ and move the contour up to $\delta$ and apply residue theorem. This involves 
obtaining rather sharp estimates on the truncated (on compact sets) resolvent near the line $\{\Re(\la)=\delta\}$. This is achieved by combining
the non-vanishing result with an a priori bound that results from a precise parametrix of the truncated resolvent.\\ 

A second result of this article is the proof of the existence of an explicit strip with infinitely many resonances.  
\begin{theo}
Let $X=\Gamma\backslash \hh^{n+1}$ be a convex co-compact hyperbolic manifold and let $\delta\in(0,n)$ be the Hausdorff
dimension of its limit set. Then for all $\varepsilon>0$, there exist infinitely many resonances in the strip 
$\{ -n\delta-\varepsilon< \Re(s)< \delta \}$. If moreover $\Gamma$ is a Schottky group, then there exist infinitely many resonances in the strip 
$\{ -\delta^2-\varepsilon< \Re(s)< \delta \}$.
\end{theo}

Note that the existence of infinitely many resonances in some strips was proved by Guillop\'e-Zworski \cite{GZw} in dimension $2$ and Perry \cite{Pe} in higher dimension, but in both cases, they did not provide any geometric information on the width of these strips. 
Our proof is based on a Selberg like trace formula and uses all previously known counting estimates for resonances. 
An interesting consequence is an explicit Omega lower bound for the remainder in \eqref{chiut}
for generic compactly supported initial data. 
\begin{cor}
For any compact set $K\subset X$, there exists a generic set $\Omega\subset C^\infty(K)$ such that for all $f_1\in\Omega,f_0=0$ and all $\eps>0$, the remainder in \eqref{chiut} is not a $\mc{O}_{L^2}(e^{-(\ndemi +n\delta+\eps)t})$ as $t\to \infty$. If $X$ is Schottky, $\mc{O}_{L^2}(e^{-(\ndemi +n\delta+\eps)t})$ can be improved to $\mc{O}_{L^2}(e^{-(\ndemi+\delta^2+\eps)t})$.
\end{cor}
The meaning of "generic" above is in the Baire category sense, i.e. it is a $G_\delta$-dense subset.
We point out that when $n=1$, all convex-cocompact surfaces are Schottky i.e. are obtained as $\Gamma\backslash \hh^2$,
where $\Gamma$ is a Schottky group. For a definition of Schottky groups in our setting we refer for example to the introduction of \cite{GLZ}. In higher dimensions, not all convex co-compact manifolds are obtained via Schottky groups. For more details and references around these questions we refer to \cite{GN}. 

The rest of the paper is organized as follows. In $\S 2$, we review and prove some necessary bounds on the resolvent in the continuation domain.
In $\S 3$ we prove the estimate on the strip with finitely many resonances. In $\S4$, we derive the asymptotics by using contour deformation and the key bounds of $\S 2$. We also show how to relate $\S 3$ to an Omega lower bound of the remainder. The section $\S 5$ is devoted to the analysis of the special cases $\delta \in  \ndemi -\nn$ in terms of the conformal theory of the infinity.

\bigskip
\noindent
\textbf{Acknowledgement}. Both authors are supported by ANR grant JC05-52556. C.G ackowledges 
support of NSF grant DMS0500788, ANR grant JC0546063 and thanks the Math department of ANU (Canberra) where part of this work was done.

\section{Resolvent}
We start in this section by analyzing the resolvent of the Laplacian for convex co-compact quotient of $\hh^{n+1}$ and 
we give some estimates of its norms.
   
\subsection{Geometric setting}   
We let $\Gamma$ be a convex co-compact group of isometries of $\hh^{n+1}$ 
with Hausdorff dimension of its limit set satisfying $0<\delta<n/2$, 
we set $X=\Gamma\backslash \hh^{n+1}$ its quotient equipped with the induced hyperbolic metric and we denote the natural projection by
\begin{equation}\label{pigamma}
\pi_\Gamma: \hh^{n+1}\to X=\Gamma\backslash \hh^{n+1}, \quad \bar{\pi}_\Gamma:
\Omega\to \pl\bar{X}=\Gamma\backslash\Omega.
\end{equation}  
By assumption on the group $\Gamma$, for any element $h\in \Gamma$
there exists $\alpha \in {\rm Isom}(\hh^{n+1})$ such that for all
$(x,y)\in \hh^{n+1}=\rr^n\times \rr_+$,
$$\alpha^{-1}\circ h \circ \alpha(x,y)=e^{l(\gamma)}(O_\gamma(x),y),$$
where $O_\gamma \in SO_n(\rr), l(\gamma)>0$. We will denote by
$\alpha_1(\gamma),\ldots,\alpha_n(\gamma)$ the eigenvalues of
$O_\gamma$, and we set
\begin{equation}\label{Ggamma}
G_\gamma(k)=\det \left(I-e^{-kl(\gamma)}O_\gamma^k \right)=
\prod_{i=1}^n \left(1-e^{-kl(\gamma)}\alpha_i(\gamma)^k \right).
\end{equation}
The Selberg zeta function of the group is defined by 
\[Z(\la)=\exp\left(-\sum_{\gamma}\sum_{m=1}^{\infty}\frac{1}{m}\frac{e^{-\la ml(\gamma)}}{G_\gamma(m)}\right),\]
the sum converges for $\Re(\la)>\delta$ and admits a meromorphic extension to $\la\in\cc$ by results
of Fried \cite{Fr} and Patterson-Perry \cite{PP}. 

\subsection{Extension of resolvent, resonances and zeros of Zeta}

The spectrum of the Laplacian $\Delta_X$ on $X$ is a half line of absolutely continuous spectrum 
$[n^2/4,\infty)$, and if 
we take for the resolvent of the Laplacian the spectral parameter $\la(n-\la)$
\[R(\la):=(\Delta_X-\la(n-\la))^{-1},\]
this is a bounded operator on $L^2(X)$ if $\Re(\la)>n/2$.
It is shown by Mazzeo-Melrose \cite{MM} and Guillop\'e-Zworski \cite{GZ2} that $R(\la)$ extends meromorphically 
in $\cc$ as continuous operators $R(\la):L^2_{\rm comp}(X)\to L^2_{\rm loc}(X)$, with poles of finite multiplicity, i.e.
the rank of the polar part in the Laurent expansion of $R(\la)$ at a pole is finite.
The poles are called \emph{resonances} of $\Delta_X$, they form the discrete set $\mc{R}$ included in $\Re(\la)<n/2$, 
where each resonance $s\in\mc{R}$ is repeated with the mutiplicity 
\[m_s:=\rang (\textrm{Res}_{\la=s}R(\la)).\]
A corollary of the analysis of divisors of $Z(\la)$ by Patterson-Perry \cite{PP} 
and Bunke-Olbrich \cite{BO} is the   
\begin{prop}[{\bf Patterson-Perry, Bunke-Olbrich}]\label{patper}
Let $s\in \cc\setminus (-\nn_0\cup(n/2-\nn))$, then $Z(\la)$ is holomorphic at $s$, 
and $s$ is a zero of $Z(\la)$ if and only if $s$ is
a resonance of $\Delta_X$. Moreover its order as zero of $Z(\la)$ coincide with the multiplicity $m_s$
of $s$ as a resonance.
\end{prop}

\subsection{Estimates on the resolvent $R(\la)$ in the non-physical sheet}
The series $P_\la(m,m')$ defined in (\ref{poincareseries}) converges absolutely in $\Re(\la)>\delta$, is a holomorphic function of $\la$ there, with local uniform bounds in $m,m'$, which clearly gives 
\[\forall \eps>0, \exists C_\eps(m,m')>0, \forall \la \textrm{ with }\Re(\la)\in[\delta+\eps,n], \quad
 |P_{\la}(m,m')|\leq C_{\eps,m,m'}\] 
and $C_{\eps,m,m'}$ is locally uniform in $m,m'$. 
We show the 
\begin{prop}\label{borneres}
With previous assumptions, there exists $\eps>0$ and a holomorphic family in $\{\Re(\la)>\delta-\eps\}$ 
of continuous operator $K(\la):L_{{\rm comp}}^2(X)\to L_{{\rm loc}}^2(X)$
such that the resolvent satisfies in $\Re(\la)>\delta$
\[R(\la)=\frac{(2\pi)^{-\ndemi}\Gamma(\la)}{\Gamma(\la-\ndemi)}P(\la)+K(\la)\]
where $P(\la)$ is the operator with Schwartz kernel $P_\la(m,m')$.
Moreover there exists $M>0$ such that for any $\chi_1,\chi_2\in C_0^\infty(X)$, there is a $C>0$ such that
\[||\chi_1K(\la)\chi_2||_{\mc{L}(L^2(X))}\leq C(|\la|+1)^M, \quad \Re(\la)>\delta-\eps\]
\end{prop}
\textsl{Proof}: we choose a fundamental domain $\mc{F}$ for $\Gamma$ with a finite 
number of sides paired by elements of $\Gamma$. 
By standard arguments of automorphic functions, the resolvent 
kernel $R(\la;m,m')$ for $m,m'\in\mc{F}$ is the average 
\[R(\la;m,m')=\sum_{\gamma\in\Gamma}G(\la;m,\gamma m')=\sum_{\gamma\in \Gamma}\sigma(d_h(m,\gamma m'))^{\la}k_\la(\sigma(d_h(m,\gamma m')))\]
\[\sigma(d):=(\cosh d)^{-1}=2e^{-d} (1+e^{-2d})^{-1}\]
where $G(\la;m,m')$ is the Green kernel of the Laplacian on $\hh^{n+1}$ and $k_\la\in C^\infty([0,1))$ is the 
hypergeometric function defined for $\Re(\la)>\frac{n-1}{2}$
\[k_\la(\sigma):=\frac{2^{\frac{3-n}{2}}\pi^{\frac{n+1}{2}}\Gamma(\la)}{\Gamma(\la-\frac{n+1}{2}+1)}
\int_0^1(2t(1-t))^{\la-\frac{n+1}{2}}(1+\sigma(1-2t))^{-\la}dt\]
which extends meromorphically to $\cc$ and whose Taylor expansion at order $2N$ can be written 
\[k_\la(\sigma)=2^{-\la-1}\sum_{j=0}^N\alpha_j(\la)\Big(\frac{\sigma}{2}\Big)^{2j}+k_\la^N(\sigma), \quad 
\alpha_j(\la):=\frac{\pi^{-\ndemi}\Gamma(\la+2j)}{\Gamma(\la-\ndemi+1)\Gamma(j+1)}\]
with $k^N_\la\in C^{\infty}([0,1))$ and the estimate for any $\eps_0>0$ 
\begin{equation}\label{estimklan}
|k_\la^N(\sigma)|\leq \sigma^{2N+2}C^N(|\la|+1)^{CN}, \quad \sigma\in[0,1-\eps_0), \quad \Re(\la)>\ndemi-N
\end{equation} 
for some $C>0$ depending only on $\eps_0$, see for instance \cite[Lem. B.1]{GTh}. 
Extracting the first term with $\alpha_0$ in $k_\la$, we can then decompose 
\[R(\la;m,m')=\frac{\pi^{\ndemi}\Gamma(\la)}{2\Gamma(\la-\ndemi+1)}\Big(
\sum_{\gamma\in\Gamma}e^{-\la d_h}+\sum_{\gamma\in\Gamma}e^{-(\la+1) d_h}
f_\la(e^{-d_h})\Big)+\sum_{\gamma\in\Gamma}\sigma(d_h)^{\la}k^0_\la(\sigma(d_h))
\]
\[f_\la(x):=\frac{(1+x^2)^{-\la}-1}{x},\]
and where $d_h$ means $d_h(m,\gamma m')$ here. Thus to prove the Proposition, we 
have to analyze the term $K(\la):=2^{-1}\alpha_0(\la)K_1(\la)+K_2(\la)$ with 
\[K_1(\la):=
\sum_{\gamma\in\Gamma}e^{-(\la+1) d_h}
f_\la(e^{-d_h}), \quad K_2(\la):=\sum_{\gamma\in\Gamma}\sigma(d_h)^{\la}k^0_\la(\sigma(d_h))\]
The first term $K_1$ is easy to deal with since $|f_\la(x)|\leq C(|\la|+1)$ for $x\in[0,1]$, thus we can use
the fact that $P_{\la+1}(m,m')$ converges absolutely in $\Re(\la)>\delta-1$, is holomorphic there, and is 
locally uniformly bounded in $(m,m')$ thus
\[|\alpha_0(\la)\chi_1(m)\chi_2(m')K_1(\la)|\leq C(|\la|+1)^{\ndemi+1}\]
the same bound holds for the operator in $\mc{L}(L^2(X))$ with Schwartz kernel 
$\chi_1(m)\chi_2(m)F_1(\la)$. Note that $\alpha_0(\la)$ has no pole in $\Re(\la)>0$, thus no pole
in $\Re(\la)>\delta/2>0$.\\

For $K_2(\la)$ we can decompose it as follows: for $m\in\textrm{Supp}(\chi_1)$, $m'\in\textrm{Supp}(\chi_2)$
(which are compact in $\mc{F}$), for $\eps_0>0$ fixed there is only a finite number of elements
$\Gamma_0=\{\gamma_0,\dots,\gamma_L\in\Gamma\}$ such that $d_h(m,\gamma m')>\eps_0$ for any 
$\gamma\notin\Gamma_0$ and any $m,m'\in\mc{F}$, this is because the group 
acts properly discontinuously on $\hh^{n+1}$.
Thus we split the sum in $K_2(\la)$ into 
\begin{equation}\label{k2}
K_2(\la)=\sum_{\gamma\in\Gamma_0}\sigma(d_h)^{\la}k^0_\la(\sigma(d_h))+\sum_{\gamma\notin \Gamma_0}
\sigma(d_h)^{\la}k^0_\la(\sigma(d_h)).\end{equation}
We first observe that the second term is a convergent series, holomorphic in $\la$, 
for $\Re(\la)>\delta-1$ and locally uniformly bounded in $(m,m')$. Indeed it is easily seen to be bounded by 
\begin{equation}\label{k22}
CN(|\la|+1)\sum_{j=1}^N|\alpha_j(\la)|P_{\Re(\la)+2j}(m,m')
+C^N(|\la|+1)^{CN}P_{\Re(\la)+2N+1}(m,m')\end{equation} 
by assumption on $\Gamma_0$ and using (\ref{estimklan}), $C$ depending on $\eps_0$ only. 
Moreover since $\alpha_j(\la)$ is polynomially bounded by $C(|\la|+1)^{2j}$
we have a polynomial bound for (\ref{k22}) of degree depending on $N$.
The first term in (\ref{k2}) has a finite sum thus it suffices to estimate each term,
but because of the usual conormal singularity of the resolvent at the diagonal, it explodes as $d_h(m,m')\to 0$.
We want to use Schur's lemma for instance, so we have to bound 
\[\sup_{m\in\mc{F}}\int_{\mc{F}}|\chi_1(m)\chi_2(m')K_2(\la;m,m')|dm'_{\hh^{n+1}},\quad
\sup_{m'\in\mc{F}}\int_{\mc{F}}|\chi_1(m)\chi_2(m')K_2(\la;m,m')|dm_{\hh^{n+1}}.
\]
First we recall that $\hh^{n+1}=(0,\infty)_x\x\rr^n_y$ has a Lie group structure with product 
\[(x,y).(x',y')=(xx',y+xy'), \quad (x,y)^{-1}=(\frac{1}{x},-\frac{y}{x})\] 
and neutral element $e:=(1,0)$. Then if $(u,v):=(x',y')^{-1}.(x,y)=(x/x',(y-y')/x')$ we get
\begin{equation}\label{coshd}
(\cosh(d_h(x,y;x',y')))^{-1}=\frac{2xx'}{x^2+{x'}^2+|y-y'|^2}=\frac{2u}{1+u+|v|^2}=(\cosh(d_h(u,v;1,0)))^{-1}.
\end{equation}
Moreover the diffeomorphism $(u,v)\to m'=m.(u,v)^{-1}$ on $\hh^{n+1}$ pulls the hyperbolic measure $dm'_{\hh^{n+1}}={x'}^{-n-1}dx'dy'$ back into the right invariant measure $u^{-1}dudv$ for the group action. This is to say that we have to bound 
\begin{equation}\label{sup1}
\sup_{m\in\mc{F}}\int_{\mc{F}^{-1}.m}|\chi_1(m)\chi_2(m.(u,v)^{-1})K_2(\la;m,m.(u,v)^{-1})|\frac{dudv}{u}\end{equation}
where $\mc{F}^{-1}.m:=\{{m'}^{-1}.m; m'\in\mc{F}\}$ and similarly
\begin{equation}\label{sup2}
\sup_{m'\in\mc{F}}\int_{\mc{F}^{-1}.m'}|\chi_1(m'.(u,v)^{-1})\chi_2(m')K_2(\la;m'.(u,v)^{-1},m')|\frac{dudv}{u}.
\end{equation}  
Because $m,m'$ are in compact sets, the estimate (\ref{k22}) with $N=n$ gives a polynomial bounds in $\la$ in 
$\{\Re(\la)>\delta-\eps\}$ 
for the terms coming from $\gamma\notin\Gamma_0$.
To deal with the term of (\ref{k2}) containing elements $\gamma\in\Gamma_0$, we use Lemma B.1 of \cite{GTh} which 
proves that for any compact $K$ of $\hh^{n+1}$, there exists a constant $C_K$ such that
\begin{equation}\label{borneresl1}
\int_{K}|G(\la;(u,v),e)|\frac{dudv}{u}\leq \frac{C_K^N(|\la|+1)^{n-1}}{\textrm{dist}(\la,-\nn_0)}, \quad \Re(\la)>\ndemi-N.\end{equation}
Now to bound (\ref{sup1}) with $K_2(\la,\bullet,\bullet)$ replaced by $\sigma(d_h(\bullet,\gamma\bullet))^\la k'_\la(\sigma(d_h(\bullet,\gamma\bullet)))$ we note that before we did our change of variable in (\ref{sup1}), 
we can make the change of variable $m'\to\gamma^{-1}m'$ which amounts to bound 
\[\sup_{m\in\mc{F}}\int_{(\gamma\mc{F})^{-1}.m}\Big|\chi_1(m)\chi_2(\gamma^{-1}m.(u,v)^{-1})\Big(G(\la;(u,v),e)-
2^{-\la-1}\alpha_0(\la)\sigma^\la(d_h((u,v),e))\Big)\Big|\frac{dudv}{u}\]
where we used (\ref{coshd}).
But again, since $\chi_1,\chi_2$ have compact support, we get a polynomial bound in $\la$ using (\ref{borneresl1})
and a trivial polynomial bound for $k_\la(0)$. 
The term (\ref{sup2}) can be dealt with similarly and we finally deduce that for some $M$,
\[||\chi_1 K_2(\la)\chi_2||_{\mc{L}(L^2(X))}\leq C(|\la|+1)^M, \quad \Re(\la)>\delta-\eps\]
and the Proposition is proved.
\qed\\ 

This clearly shows that the resolvent extends to $\{\Re(\la)>\delta\}$ analytically. Actually, 
Patterson \cite{Pa} (see also \cite[Prop 1.1]{P}) showed the following.
\begin{prop}[{\bf Patterson}]\label{pa}
The family of operators $\Gamma(\la-n/2+1)R(\la)$ is holomorphic in $\{\Re(\la)>\delta\}$,
has no pole on $\{\Re(\la)=\delta, \la\not=\delta\}$ and has a pole of order $1$ at $\la=\delta$
with rank $1$ residue given by
\[{\rm Res}_{\la=\delta} \Gamma(\la-n/2+1)R(\la)=A_X u_\delta\otimes u_\delta\]
where $A_X\not=0$ is some constant depending on $\Gamma$ and $u_\delta$
is the Patterson generalized eigenfunction defined by
\begin{equation}\label{udelta}
\pi_\Gamma^*u_\delta(m)=\int_{\pl_\infty\hh^{n+1}}\Big(\mc{P}(m,y)\Big)^\delta d\mu_{\Gamma}(y)
\end{equation} 
$\mc{P}$ being the Poison kernel of $\hh^{n+1}$ and $d\mu_{\Gamma}$ the 
Patterson-Sullivan measure associated to $\Gamma$ on the sphere $\pl_{\infty}\hh^{n+1}=\rr^n\cup\{\infty\}
= S^n$.
\end{prop}

We can but notice that $\delta\in n/2-\nn$ is a special case since the resolvent becomes holomorphic at $\la=\delta$.
We postpone the analysis of this phenomenon to section $\S 5$.\\

A rough exponential estimate in the non-physical sheet also holds using determinants methods
 (used for instance in \cite{GZ}). 
\begin{lem}\label{exponentres}
For $\chi_1,\chi_2\in C_0^\infty(X)$, $j\in\nn_0$, and $\eta>0$ there is $C>0$ such that for $|\la|\leq N/16$ 
and $\textrm{dist}(\la,\mc{R}\}>\eta$, 
\[\quad ||\pl_\la^j\chi_1R(\la)\chi_2||_{\mc{L}(L^2(X))}\leq e^{C(N+1)^{n+3}},\]
\end{lem}
\textsl{Proof}: we apply the idea of \cite[Lem. 3.6]{GZ} with the parametrix construction
of $R(\la)$ written in \cite{GZ2}. Let $x$ be a boundary defining function 
of $\pl\bar{X}$ in $\bar{X}$, which can be considered as a weight to define Hilbert spaces
$x^\alpha L^2(X)$, for any $\alpha\in\rr$. 
For any large $N>0$ that we suppose in $2\nn$ for convenience, 
Guillop\'e and Zworski \cite{GZ2} construct operators 
\[P_N(\la,\la_0):x^NL^2(X)\to x^{-N}L^2(X), \quad K_N(\la,\la_0):x^NL^2(X)\to X^NL^2(X)\] 
meromorphic with finite multiplicity in $O_N:=\{\Re(\la)>(n-N)/2\}$, whose poles are situated
at $-\nn_0$, and such that
\[(\Delta_X-\la(n-\la))P_N(\la,\la_0)=1+K_N(\la,\la_0)\]
with $\la_0$ large depending on $N$, take for instance $\la_0=n/2+N/8$. 
Moreover $K_N(\la,\la_0)$ is compact with characteristic values satisfying in 
$O_{N,\eta}:=O_N\cap\{\textrm{dist}(\la,-\nn_0)>\eta\}$ 
\begin{equation}\label{muj}
\mu_j(K_N(\la,\la_0))\leq C(1+|\la-\la_0|)j^{-\frac{1}{n}}+\left\{\begin{array}{ll}
e^{CN} & \textrm{ if } j\leq CN^{n+1}\\
e^{-N/C}j^2  & \textrm{ if } j\geq CN^{n+1}
\end{array}\right.
\end{equation}
for some $0<\eta<1/4$ and $C>0$ independent of $\la,N$.
They also have $||K_N(\la_0,\la_0)||\leq 1/2$ in $\mc{L}(x^NL^2(X))$,
thus by Fredholm theorem 
\[R(\la)=P_N(\la,\la_0)(1+K_N(\la,\la_0))^{-1}: x^NL^2(X)\to x^{-N}L^2(X)\]
is meromorphic with poles of finite multiplicity in $O_N$.
By standard method as in \cite[Lem. 3.6]{GZ} we define 
\[d_N(\la):=\det (1+K_N(\la,\la_0)^{n+2})\] 
which exists in view of (\ref{muj}), and we have the rough bound 
\begin{equation}\label{estimekn}
||(1+K_N(\la,\la_0))^{-1}||_{\mc{L}(x^NL^2(X))}\leq \frac{\det(1+|K_N(\la,\la_0)|^{n+2})}{|d_N(\la)|}
\end{equation}
in $O_{N,\eta}$ and where $|A|:=(A^*A)^{\demi}$ for $A$ compact.
The term in the numerator is easily shown to be bounded by $\exp(C(N+1)^{n+2})$ in $O_{N,\eta}$ from 
(\ref{muj}), actually this is written in \cite[Lem. 5.2]{GZ2}. It remains to have a lower bound of $|d_N(\la)|$.
In Lemma 3.6 of \cite{GZ}, they use the minimum modulus theorem to obtain lower bound of a function 
using its upper bound, but this means that the function has to be analytic in $\cc$. 
Here there is a substitute which is Cartan's estimate \cite[Th. I.11]{Le}. 
We first need to multiply $d_N(\la)$ by a holomorphic function $J_N(\la)$
with zeros of sufficient multiplicity at $\{-k;k=0,\dots,N/2\}$ 
to make $J_N(\la)d_N(\la)$ holomorphic in $O_N$, for instance the polynomial
\[J_N(\la):=\prod_{k=0}^{N/2}(\la-k)^{CN^{n+2}}\]
for some large integer $C>0$ suffices in view of the order ($\leq C{N^{n+2}}$) of each $-k$ 
as a pole of $d_N(\la)$ proved in \cite[Lem. A.1]{GZ2}. Then clearly $f_N(\la):=J_N(\la+\la_0)d_N(\la+\la_0)/(J_N(\la_0)d_N(\la_0))$
is holomorphic in $\{|\la|\leq N/4\}$ and satisfies in this disk 
\[|f_N(\la)|\leq e^{C(N+1)^{n+3}},  \quad f_N(0)=1,\]
where we used the maximum principle in disks around each $-k$ to estimate the norm there.
Thus we may apply Cartan's estimate for this function in $|\la|<N/4$: for all $\alpha>0$ small enough 
there exists $C_\alpha>0$ such that, outside a family of disks
the sum of whose radii is bounded by $\alpha N$
\[\log|f_N(\la)|>-C_\alpha\log \Big(\sup_{|\la|\leq N/4}|f_N(\la)|\Big)\] 
and $|\la|\leq N/4$.
Fixing $\alpha$ sufficiently small, there exists $\beta_N\in(3/4,1)$ so that
\[|d_N(\la)|>e^{-C(N+1)^{n+3}} \textrm{ for }|\la-\la_0|=\beta_N \frac{N}{4}.\]
Note that we can also choose $\beta_N$ so that $\textrm{dist}(\beta_NN/4,\nn)>\eta$ for some small $\eta$ 
uniform with respect to $N$.
Thus the same bound holds for $||(1+K_N(\la,\la_0))^{-1}||_{\mc{L}(x^NL^2(X))}$ using (\ref{estimekn}).
Now we need a bound for $P_N(\la,\la_0)$ and
it suffices to get back to its definition in the proof of Proposition 3.2 of \cite{GZ2}: it involves operators 
of the form $\iota^*\varphi R_{\hh^{n+1}}(\la)\psi\iota_*$ for some cut-off functions $\psi,\varphi \in C^\infty
(\hh^{n+1})$ and isometry  
\[\iota: U\subset X\to \{(x,y)\in (0,\infty)\x \rr^n; x^2+|y|^2<1\}\subset \hh^{n+1},\] 
and operators whose norm is explicitely bounded in \cite[Sect. 4]{GZ2} by $e^{C(N+1)}$ in $O_{N,\eta}$.
The appendix B of \cite{GTh} gives an estimate of the same form for $||\varphi R_{\hh^{n+1}}(\la)\psi||$ 
as an operator in $\mc{L}(x^NL^2(X), x^{-N}L^2(X))$ for $\la\in O_{N,\eta}$ 
(this is actually a direct consequence of (\ref{borneresl1}) and (\ref{estimklan})) thus we have the bound 
\[||R(\la)||_{\mc{L}(x^NL^2(X), x^{-N}L^2(X))}\leq e^{C(N+1)^3}\]  
in $\{|\la-\la_0|=\beta_NN/4\}$. 
Let $\mc{R}_N$ be the set of poles of $R(\la)$ in $O_{N}$, each pole being repeated according to its order; 
$\mc{R}_N$ has at most $CN^{n+2}$ elements so we may multiply
$R(\la)$ by 
\[F_N(\la):=\prod_{s\in \mc{R}_N}E(\la/s,n+2)\]
where $E(z,p):=(1-z)\exp(z+\dots+p^{-1}z^p)$ is the Weierstrass elementary function.
It is rather easy to check that for all $\eps>0$ small, we have the bounds
\begin{equation}\label{estimateF_N}
 e^{C_\eps (N+1)^{n+3}}\geq |F_N(\la)|\geq e^{-C_\eps(N+1)^{(n+3)}}
 \end{equation}
for some $C_\eps$ and for all $\la\in O_{N}$ such that $\textrm{dist}(\la,\mc{R})>\eps$.
Thus $R(\la)F_N(\la)$ is holomorphic in $\{|\la-\la_0|\leq \beta_NN/4\}$ and we can use the maximum
principle which gives a upper bound $||F_N(\la)R(\la)||_{\mc{L}(x^NL^2, x^{-N}L^2)}
\leq \exp(C_\eps(N+1)^{n+3})$ in 
$\{|\la-\la_0|\leq \beta_NN/4\}$. We get our conclusion using \eqref{estimateF_N}, the fact that $\chi_i$
is bounded by $e^{CN}$ as an operator from $L^2$ to $x^NL^2$, and the Cauchy formula for the case $j>0$ 
(estimates of the derivatives with respect to $\la$). 
\qed\\

\textsl{Remark}: Notice that similar estimates are obtained independently by Borthwick \cite{Bor}.\\

In the case of surfaces the second author \cite{N} used the powerful estimates 
developped by Dolgopyat \cite{Do} to prove that the Selberg zeta function $Z(\la)$ is analytic and non-vanishing in 
$\{\Re(\la)>\delta-\eps, \la\not=\delta\}$ for some $\eps>0$.
In higher dimension, the same result holds, as was shown recently by Stoyanov \cite{S}. 
\begin{theo}[\bf{Naud, Stoyanov}]\label{naudstoy}
There exists $\eps>0$ such that the Selberg zeta function $Z(\la)$ is holomorphic and non-vanishing in 
$\{\la\in\cc; \Re(\la)>\delta-\eps, \la\not=\delta\}$. 
\end{theo}

Using Proposition \ref{patper}, this implies that the resolvent $R(\la)$ is holomorphic in 
a similar set (possibly by taking $\eps>0$ smaller). 
Then an easy consequence of the maximum principle as in \cite{TZ,BP} together with a
rough exponential bound for the resolvent allows to get a polynomial 
bound for $||\chi_1R(\la)\chi_2||$ on the $\{\Re(\la)=\delta; \la\not=\delta\}$. 
 
\begin{cor}\label{extres}
There is $\eps>0$ such that the resolvent $R(\la)$ is meromorphic in $\Re(\la)>\delta-\eps$ 
with only possible pole the simple pole $\la=\delta$, the residue of which is given by
\[\textrm{Res}_{\la=\delta}R(\la)=\frac{A_X}{\Gamma(\ndemi-\delta+1)}u_\delta\otimes u_\delta\]
where $u_\delta$ is the Patterson generalized eigenfunction of (\ref{udelta}), $A_X\not=0$ a constant.
Moreover for all $j\in\nn_0$, $\chi_1,\chi_2\in C_0^\infty(X)$, there exists $L\in\nn, C>0$ such that for $|\la-\delta|>1$ 
\[||\pl^j_\la\chi_1R(\la)\chi_2||_{\mc{L}(L^2(X))}\leq C(|\la|+1)^{L}
\textrm{ in }\{\Re(\la)\geq \delta\}\]  
\end{cor}
\textsl{Proof}: This  is a consequence of 
Proposition \ref{borneres}, Proposition \ref{pa}, Theorem \ref{naudstoy}   
and the maximum principle as in \cite[Prop. 1]{BP}.
First we remark from Proposition \ref{borneres} and Proposition \ref{pa} that $P_\la$ has a first order pole
with rank one residue at $\la=\delta$ and, since $|P_\la(m,m')|\leq |P_{\Re(\la)}(m,m')|$, 
we have the estimate
\[||\chi_1R(\la)\chi_2||_{\mc{L}(L^2(X))}\leq |\Re(\la)-\delta|^{-1}C(|\la|+1)^{M}\]
for $\Re(\la)\in (\delta,n/2)$. This implies by the Cauchy formula that 
\[||\pl_\la^j\chi_1R(\la)\chi_2||_{\mc{L}(L^2(X))}\leq |\Re(\la)-\delta|^{-1-j}C(|\la|+1)^{M}.\]
Let $A>0$, and $\varphi,\psi\in L^2(X)$, we can apply the maximum principle to the function 
\[f(\la)=e^{iA(-i(\la-\delta))^{n+4}} \cjg \pl^j_\la\chi_1R(\la)\chi_2\varphi,\psi\cjd\]
which is holomorphic in the domain $\Lambda$ bounded by the curves 
\[\Lambda_+:=\{\delta+u^{-n-3}+iu; u>1\},
\quad \Lambda_-:=\{\delta-\eps+iu;u>1\}, \quad \Lambda_0:=\{i+u; \delta-\eps<u<\delta+1\}.\]
Then it is easy to check as in \cite[Prop. 1]{BP} that by choosing $A>0$ large enough
\[ |f(\la)|<C(|\la|+1)^L||\varphi||_{L^2}||\psi||_{L^2}\] 
in $\Lambda$ for some $L$ depending only on $M$.
In particular, applying the same method in the symmetric domain $\bar{\Lambda}:=\{\bar{\la};\la\in\Lambda\}$, we obtain the polynomial bound $||\pl_\la^j\chi_1R(\la)\chi_2||\leq C(|\la|+1)^L$ on $\{\Re(\la)=\delta,|\Im(\la)|>1\}$.
\qed

\section{Width of the strip with finitely many resonances}

As stated in Theorem \ref{naudstoy}, we know that there exists a strip $\{\delta-\eps<\Re(\la)<\delta\}$
with no resonance for $\Delta_g$, or equivalently no zero for Selberg zeta function. However the proof of this 
result does not provide any effective estimate on the width of this strip (i.e. on $\eps$ above).
More generally it is of interest to know the following 
\[\rho_\Gamma:=\inf \Big\{ s\in\rr; Z(\la) \textrm{ has at most finitely many zero in }\{\Re(\la)>s\}\Big\}\]
or equivalently
\[\rho_\Gamma = \inf \Big\{ s\in\rr; R(\la) \textrm{ has at most finitely many poles in }\{\Re(\la)>s\}\Big\}.\]
In this work, we give a lower bound for $\rho_\Gamma$:   
\begin{theo}\label{naudstriptease}
Let $X=\Gamma\backslash \hh^{n+1}$ be a convex co-compact hyperbolic manifold and let $\delta\in(0,n)$ be the Hausdorff
dimension of its limit set. Then for all $\varepsilon>0$, there exist infinitely many resonances in the strip 
$\{ -n\delta-\varepsilon< \Re(s)< \delta \}$. If moreover $\Gamma$ is a Schottky group, then there exist infinitely many resonances in the strip 
$\{ -\delta^2-\varepsilon< \Re(s)< \delta \}$.
\end{theo}
\noindent\textsl{Remark}: In particular, we have $\rho_\Gamma\geq-\delta n$ in general and
$\rho_{\Gamma}\geq -\delta^2$ for Schottky manifolds.
The limit case $\delta\rightarrow 0$ may be viewed as
a cyclic elementary group $\Gamma_0$, and resonances of the Laplace operator on 
$\Gamma_0\backslash \H$ are given explicitely \cite[Appendix]{GZ1}, they form a lattice 
$\{-k+i\alpha \ell; k\in\nn_0,\ell\in\zz\}$ for some $\alpha\in\rr$, 
in particular there are infinitely many resonances on the vertical line
$\{ \Re(s)=0 \}$. This heuristic consideration suggests that for small values of $\delta$, 
our result is rather sharp.\\

\textsl{Proof}: 
The proof is based on the trace formula of \cite{GN} and estimates 
on the distribution of resonances due to Patterson-Perry \cite{PP}, Guillop\'e-Lin-Zworski \cite{GLZ} (see also Zworski \cite{Z} for dimension $2$).
To make some computations clearer (Fourier transforms), we will use the spectral parameter $z$ with
$\la=\ndemi+iz$ and $\Im z >0$ in the non-physical half-plane. 
We set $\beta:=\delta$ if $X$ is Schottky or $n+1=2$, while $\beta:=n$ if $n+1>2$ and $X$ not Schottky.
We proceed by contradiction and assume that there is $\rho=n/2+\beta\delta+\varepsilon$ for some $\varepsilon>0$ such that
there are at most finitely many resonances in $\Im(z)<\rho$. 
Let us first recall the trace formula of \cite{GN}: as distributions of $t\in\rr\setminus\{0\}$, we have the identity 
\begin{equation}\label{trace}
\begin{gathered}
\demi\Big(\sum_{\ndemi+iz\in\mc{R}}e^{iz|t|}+\sum_{k\in\nn}d_ke^{-k|t|}\Big)
=\sum_{\gamma\in\mc{P}}\sum_{m=1}^{\infty}\frac{\ell(\gamma)e^{-\ndemi m\ell(\gamma)}
}{2G_\gamma(m)}\delta(|t|-m\ell(\gamma))+\frac{\chi(\bar{X})\cosh\frac{t}{2}}
{(2\sinh\frac{|t|}{2})^{n+1}},
\end{gathered}
\end{equation}
where ${\mc P}$ denotes the set of primitive
closed geodesics on $X=\Gamma\backslash \hh^{n+1}$, $\ell(\gamma)$ stands for the length of $\gamma\in\mc{P}$, $G_\gamma(m)$ is defined in \eqref{Ggamma},
$d_k:=\dim\ker P_k$ if $P_k$ is the k-th GJMS conformal Laplacian on the conformal boundary $\pl\bar{X}$, $\mc{R}$ is the set of resonances of $\Delta_X$ counted with multiplicity and $\chi(\bar{X})$ denotes the Euler characteristic of $\bar{X}$.
\noindent
Next we choose $\varphi_0\in C_0^\infty(\R)$ a positive weight supported on $[-1,+1]$
with $\varphi_0(0)=1$ and $0\leq \varphi_0 \leq 1$. We set 
$$\varphi_{\alpha,d}(t)=\varphi_0\left(\frac{t-d}{\alpha}\right),$$
where $d$ will be a large positive number and $\alpha>0$ will be small when compared to 
$d$ (typically $\alpha=e^{-\mu d}$). 
Pluging it into the trace formula \eqref{trace} and assuming that $d$ coincides with 
a large length of a closed geodesic, we get that for $d$ large enough,
$$\sum_{\gamma,m}\frac{\ell(\gamma)e^{-\ndemi m \ell(\gamma)}}{2G_{\gamma}(m)}\varphi_{\alpha,d}(ml(\gamma))
\geq Ce^{-\ndemi d},$$
with a constant $C>0$, whereas the other term can be estimated by 
\[\alpha\chi(\bar{X})\int_{-1}^1\varphi(t)\frac{\cosh((d+t\alpha)/2)}
{(2\sinh(|d+t\alpha|/2))^{n+1}}dt=\mc{O}(\alpha)e^{-\ndemi d}.\] 
The key part of the proof is to estimate carefully the spectral side
of the formula, i.e. we must examinate 
$$\sum_{\ndemi+iz \in {\mc R} } \widehat{\varphi}_{\alpha,d}(-z)+\sum_{\substack{\ndemi+iz=-k\\
k\in\nn_0}} d_k\widehat{\varphi}_{\alpha,d}(-z),$$
where $\widehat{\varphi}$ is the usual Fourier transform. Standard formulas on Fourier 
transform on the Schwartz space show that for all integer $M>0$, there exists a constant 
$C_M>0$ such that 
\begin{equation}
\label{est1}
\left|\widehat{\varphi}_{\alpha,d}(-z) \right| \leq \alpha C_M 
\frac{e^{-d\Im(z)+\alpha|\Im(z)|}}{(1+\alpha|z|)^M}.
\end{equation}
To simplify, we denote by $\widetilde{{\mc R}}$ the set 
$\{z \in \C;\ndemi+iz \in {\mc R}\cup i\nn\}$ where 
each element $z$ is repeated with the multiplicity 
\[\left\{\begin{array}{l}
m_{n/2+iz} \textrm{ if }z\notin i\nn\\
m_{n/2-k}+d_k \textrm{ if }z=ik \textrm{ with }k\in\nn
\end{array}\right..\]
Our assumption now is that
$$\{ 0\leq \Im(z) \leq \rho\}\cap \widetilde{{\mc R}}$$
is finite for $\rho=\ndemi+\beta\delta+\varepsilon$. We set $\overline{\rho}>\rho\geq 0$. The idea is to split the sum over 
resonances as
$$\sum_{z \in \widetilde{{\mc R}_X}} \widehat{\varphi}_{\alpha,d}(-z)
=\sum_{\ndemi-\delta\leq \Im(z) \leq \rho} \widehat{\varphi}_{\alpha,d}(-z) 
+\sum_{\rho\leq \Im(z) \leq \overline{\rho}} \widehat{\varphi}_{\alpha,d}(-z) 
+\sum_{\overline{\rho}\leq \Im(z)} \widehat{\varphi}_{\alpha,d}(-z),$$  
and estimate their contributions using dimensional and fractal upper bounds.
Using (\ref{est1}) we can bound the last term (for $d$ large) by
$$\left|\sum_{\overline{\rho}\leq \Im(z)} \widehat{\varphi}_{\alpha,d}(-z)
\right| \leq C_M \alpha e^{-\overline{\rho}(d-\alpha)}
\int_{\overline{\rho}}^{+\infty} \frac{d{\mc N}(r)}{(1+\alpha r)^M},$$
where ${\mc N}(r)=\#\{z \in \widetilde{{\mc R}}; |z|\leq r\}$.
By \cite[Th. 1.10]{PP} (see also \cite[Lemma 2.3]{GN} for a discussion about the $d_k$ terms), 
we know that ${\mc N}(r)=\mc{O}(r^{n+1})$, thus we can choose $M=n+2$ and obtain, after
a Stieltjes integration by parts, the following upper bound
$$\left|\sum_{\overline{\rho}\leq \Im(z)} \widehat{\varphi}_{\alpha,d}(-z)
\right|=\mc{O}(\alpha^{-n}e^{-\overline{\rho}d} ).$$
Similarly, we have the estimate (for $d$ large and $\alpha$ small)
$$\left|\sum_{\rho\leq \Im(z) \leq \overline{\rho}} \widehat{\varphi}_{\alpha,d}(-z) \right|\leq C_M \alpha e^{-\rho(d-\alpha)} 
\int_{\rho}^{+\infty} \frac{d\widetilde{\mc N}(r)}{(1+\alpha r)^M},$$
where $\widetilde{\mc N}(r)=
\#\{z \in \widetilde{{\mc R}} \ :\ \rho\leq  \Im(z) \leq 
\overline{\rho},\ |z|\leq r\}$. This counting function is known to enjoy the 
``fractal'' upper bound $\widetilde{\mc N}(r)=\mc{O}(r^{1+\delta})$ when  $X$ is Schottky \cite{GLZ} (see also  
\cite{Z} when $n=1$), thus we can write $\widetilde{\mc N}(r)=\mc{O}(r^{1+\beta})$ where $\beta$ is defined above.
In other words, one obtains by choosing $M=n+2$, 
$$\left|\sum_{\rho\leq \Im(z) \leq \overline{\rho}} 
\widehat{\varphi}_{\alpha,d}(-z) \right|=\mc{O}(\alpha^{-\beta}e^{-\rho d} ).$$ 
Since we have assumed that $\{ 0\leq \Im(z) \leq \rho\}\cap \widetilde{{\mc R}}$
is finite, and using the fact that resonances (in the $z$ plane) have all
imaginary part greater than $\ndemi-\delta$, we also get
$$\left|\sum_{\ndemi-\delta\leq \Im(z) \leq \rho} 
\widehat{\varphi}_{\alpha,d}(-z)\right| =\mc{O}(\alpha e^{(\delta-\ndemi)d} ).$$
Gathering all estimates, we have obtained as $d \rightarrow +\infty$,
$$e^{-\ndemi d}(C+\mc{O}(\alpha))=\mc{O}(\alpha e^{(\delta-\ndemi)d} )+\mc{O}(\alpha^{-\beta}e^{-\rho d} )
+\mc{O}(\alpha^{-n}e^{-\overline{\rho}d} ),$$
where all the implied constants do not depend on $d$ and $\alpha$. If we now
set $\alpha=e^{-\mu d}$, we get a {\it contradiction} as $d\rightarrow +\infty$, provided that
$$\left \{ \begin{array}{ccc}
n\mu-\overline{\rho}&<&-\ndemi\\
\delta&<& \mu\\
\rho-\beta \mu&>& \ndemi.
\end{array} \right.$$
Set $\mu:=\delta+\varepsilon$ and
$\rho=\beta \mu+\ndemi +\varepsilon=\beta\delta+\ndemi +\varepsilon(1+\beta)$, we can then choose 
$\overline{\rho}:=n\mu+\ndemi+2\varepsilon$ which is larger than $\rho$ and we have our contradiction for all $\varepsilon>0$.\qed\\
 
The proof reveals that any precise knowledge in the 
asymptotic distribution of resonances in strips has a direct impact on resonances with small imaginary part.

\section{Wave asymptotic}
\subsection{The leading term}
Let $f,\chi\in C_0^\infty(X)$, it is sufficient to describe the large time asymptotic of the function
\[u(t):=\chi \frac{\sin(t\sqrt{\Delta_X-\frac{n^2}{4}})}{\sqrt{\Delta_X-\frac{n^2}{4}}}f\]
and $\pl_tu(t)$.
We proceed using same ideas than in \cite{CZ}. We first recall that from Stone formula the spectral measure is
\[d\Pi(v^2)=\frac{i}{2\pi}\Big(R(\ndemi+iv)-R(\ndemi-iv)\Big)dv\]
in the sense that for $h\in C^\infty([0,\infty))$ we have 
\[h\Big(\Delta_X-\frac{n^2}{4}\Big)=\int_{0}^\infty h(v)d\Pi(v^2)2vdv.\]
Since $\sin$ is odd, then it is clear that $u(t)$ can be expressed by the integral 
\begin{equation}\label{ut}
u(t)=\frac{1}{2\pi}\int_{-\infty}^\infty e^{itv}\Big(\chi R(\ndemi+iv)f-\chi R(\ndemi-iv)f\Big)dv
\end{equation}
which is actually convergent since $f\in C_0^\infty$ (this is shown below).
We want to move the contour of integration into the non-physical sheet $\{\Im(v)>0\}$ (which correponds
with $\la=n/2+iv$ to $\{\Re(\la)<n/2\}$) for the part with $e^{itv}$ and into 
the physical sheet $\{\Im(v)<0\}$ for the part with $e^{-itv}$. After setting 
\[L(v):=\Big(\chi R(\ndemi+iv)f-\chi R(\ndemi-iv)f\Big)\]
and $\eta>0$ small, we study the following integral for $\beta:=n/2-\delta$
\[I_1(R,\eta,t):=\int_{\substack{\Im(v)=\beta\\
\eta<|\Re(v)|<R}}e^{itv}L(v)dv, \quad I_2(R,t):=\int_{\substack{|\Re(v)|=R\\
0<\Im(v)<\beta}}e^{itv}L(v)dv, \]
In particular let us first show that 
\begin{lem}\label{i2}
If $|L(v)|\leq C(|v|+1)^{M}$ in $\{|\Im (v)|\leq\beta\}$ for some $C,M>0$, then 
\[\lim_{R\to\infty}I_2(R,t)=\lim_{R\to\infty} \pl_tI_2(R,t)=0\] 
\end{lem}
\textsl{Proof}: it suffices to prove inverse polynomial bounds for $L(v)$ as $|\Re(v)|\to \infty$.
Actually we can rewrite $L(v)$ using Green formula \cite{G1,P}
\begin{equation}\label{ltv}
L(v;m)=-2iv\int_X\int_{\pl\bar{X}}\chi(m)E\Big(\ndemi+iv;m,y\Big)E\Big(\ndemi-iv;m',y\Big)f(m')dy_{\pl\bar{X}}dm'_X
\end{equation}
where $E(\la;m,y)$ denotes the Eisenstein function, or equivalently the Schwartz kernel 
of the Poisson operator (see \cite{JSB}), they satisfy for all $y\in\pl\bar{X}$ 
\[\Big(\Delta_X-\frac{n^2}{4}-v^2\Big)E\Big(\ndemi+iv;\bullet,y\Big)=0.\]
Using this equation, integrating by parts $N$ times in $m'$ in (\ref{ltv}) and using the assumed polynomial bound
of $|L(v)|$ in $|\Im(v)|\leq \beta$, we get for all $N>0$ (recall $f\in C_0^\infty(X)$)
\begin{equation}\label{borneltv}
|L(v)|\leq C_N(|\Re(v)|+1)^{M-N}\end{equation} 
for some constant $C_N$. Then it suffices to take $N\gg M$ large enough and the Lemma is proved.
\qed\\

Now we get estimates in $t$ for $I_1(R,\eta,t)$.
\begin{lem}\label{i1}
If $|L(v)|\leq C(|v|+1)^{M}$ in $|\Im v|\leq\beta$ for some $C,M>0$, 
then $I_1(R,\eta,t)$ and $\pl_tI_1(R,\eta,t)$ have a limit as $R\to\infty,\eta\to 0$ and
\[\lim_{\eta\to 0}\lim_{R\to \infty}I_1(R,\eta,t)=\pi ie^{-\beta t}{\rm Res }_{v=i\beta}(L(v))+\mc{O}(e^{-\beta t}t^{-\infty}), \quad t\to \infty,\]
\[\lim_{\eta\to 0}\lim_{R\to \infty}\pl_tI_1(R,\eta,t)=-\pi\beta ie^{-\beta t}{\rm Res }_{v=i\beta}(L(v))+\mc{O}(e^{-\beta t}t^{-\infty}), \quad t\to \infty\] 
\end{lem}
\textsl{Proof}: Let us first consider $I_1(R,\eta,t)$, it can clearly be written as
\[e^{-t\beta}\int_{\eta<|u|<R}e^{itu}L(u+i\beta)du.\]
Since $L(u+i\beta)$ has a pole at $u=0$, we can write
\[L(u+i\beta)=\frac{a}{u}+h(u)\] 
for some residue $a\in\cc$ and $h(u)$ analytic on $\rr$. Set $\psi \in C_0^\infty((-1,1))$ even and equal to $1$
near $0$, then by (\ref{borneltv}) and properties of Fourier transform the integral
\[\int_{\eta<|u|<R}e^{itu}\Big((1-\psi(u))L(u+i\beta)+ \psi(u)h(u)\Big)du, \]    
converges as $R\to \infty,\eta\to 0$ to a function that is a $\mc{O}(t^{-\infty})$ when $t\to \infty$.
Now it remains to consider
\[a\int_{\eta<|u|<R}e^{itu}\psi(u)u^{-1}du=2ia\int_{\eta}^R\frac{\sin(ut)}{u}\psi(u)du\]
which clearly has a limit as $R\to\infty,\eta\to 0$, we denote by $s(t)$ this limit.
Then since $s(0)=0$ and $\psi(-u)=\psi(u)$, we have
\[\pl_t s(t)=2ia\int_{0}^\infty \psi(u)\cos(tu)du=ia\hat{\psi}(t), \quad s(t)=ia\int_0^t \hat{\psi}(\xi)d\xi=\demi ia
\int_{-t}^t\hat{\psi}(\xi)d\xi\]
and it is clear that 
\[s(t)=\lim_{t\to \infty}s(t)+\mc{O}(t^{-\infty})=\pi ia+\mc{O}(t^{-\infty}).\]
The same arguments show that 
\[\pl_ts(t)=\mc{O}(t^{-\infty})\]
and this proves the result.  
\qed\\

Now we can conclude 
\begin{theo}\label{mainth}
Let $\chi\in C_0^\infty(X)$, then the solution $u(t)$ of we wave equation \eqref{waveeq} satisfies the asymptotic
\[\chi u(t)=\frac{A_X}{\Gamma(\ndemi-\delta+1)}e^{-t(\ndemi-\delta)}\cjg u_\delta,(\delta-n/2)f_0+f_1\cjd \chi u_\delta +
\mc{O}_{L^2}(e^{-t(\ndemi-\delta)}t^{-\infty})\] 
as $t\to +\infty$, where $u_\delta$ is the Patterson generalized eigenfunction.  
\end{theo}
\textsl{Proof}: we apply the residue theorem after changing the contour in (\ref{ut}) as explained above. 
This gives for instance for $f=(0,f_1)$,
\[\int_{-R}^Re^{itv}L(v)dv=I_1(R,\eta,t)+I_2(R,t)+\int_{\substack{v=i\beta+\eta \exp(i\theta)\\
-\pi<\theta<0}}e^{itv}L(v)dv\] 
The limit of the last integral as $\eta\to 0$ is given $\pi ie^{-\beta t}\textrm{Res}_{v=i\beta} L(v)$.
It suffices to conclude by taking the limits $R\to\infty,\eta\to 0$ and using Lemmas \ref{i2} and \ref{i1}.
Then the case $f=(f_0,0)$ is dealt with similarly by differentiating in $t$ the equation above and using Lemmas
\ref{i1}, \ref{i2}.
\qed\\

We now show a lower bound in $t$ for the remainder in $u(t)$ using Theorem \ref{naudstriptease}. 
\begin{prop}\label{rem} 
Let $K\subset X$ be a compact set, then there exists a generic set $\Omega\subset C^\infty(K)$ (i.e. a countable intersection of open dense sets) such that for all $f_1\in\Omega$ and all $\varepsilon>0$, we have $r(t)\not=\mc{O}_{L^2}(e^{-(\ndemi+n\delta+\varepsilon)t})$ where
\[r(t):=\chi u(t)-\frac{A_X}{\Gamma(\ndemi-\delta+1)}e^{-t(\ndemi-\delta)}\cjg u_\delta,f\cjd \chi u_\delta\]
is the remainder in the expansion of the solution $u(t)$ of the wave equation \eqref{waveeq} with initial data $(0,f_1)$. 
The lower bound can be improved by $r(t)\not=\mc{O}_{L^2}(e^{-(\ndemi+\delta^2+\varepsilon)t})$ if $X$ is Schottky.
\end{prop}
\textsl{Proof}: Let us define $\Omega$. If $\la_0$ is a resonance, we denote by $\Pi_{\la_0}$ the polar part
in the Laurent expansion of $R(\la)$ at $\la_0$. It is a finite rank operator of the form
\[\Pi_{\la_0}= \sum_{j=1}^k(\la-\la_0)^{-j}\sum_{m=1}^{m_j(\la_0)}\varphi_{jm}\otimes \psi_{jm}\]
where $m_j(\la_0),k\in\nn$ and $\psi_{jm},\varphi_{jm}\in C^\infty(X)$, see for instance Lemma 3.1 of \cite{G}. Thus it is a continuous operator
from $C^\infty(K)$ to $C^{\infty}(X)$ and thus the kernel of $\chi \Pi_{\la_0}|_{C^\infty(K)}$ 
is a closed nowhere dense set of $C^\infty(K)$, we thus define 
$\Omega=\cap_{s\in\mc{R}}(C^\infty(K)\setminus\ker \chi\Pi_{s}|_{C^\infty(K)})$ which is a generic set of $C^\infty(K)$ (recall that $C^\infty(K)$
is a Frechet space by compactness of $K$).
The idea now is to use the existence of a resonance, say $\la_0$, in the strip
$\{-n\delta+\varepsilon>\Re(\la)>\delta\}$ proved in Theorem \ref{naudstriptease} 
and the formula (for $\Re(\la)>\delta$)
\[\chi R(\la)f=\int_{0}^\infty e^{t(\ndemi-\la)}\chi u(t)dt.\]
Indeed, if $r(t)=\mc{O}(e^{-t(\ndemi+n\delta+\varepsilon)})$, the integral $\int_{0}^\infty e^{t(\ndemi-\la)}r(t)dt$
converges for $\Re(\la)>-n\delta-\varepsilon$, and so it provides a holomorphic continuation of $\chi R(\la)f$ in 
$\la$ there. Now a straightforward computation combined with Corollary \ref{extres} 
shows that for $\Re(\la)>\delta$
\[\int_{0}^\infty e^{(\ndemi-\la)t}r(t)dt=\chi R(\la)f-(\la-\delta)^{-1}\chi\textrm{Res}_{\la=\delta}R(\la)f.\]
This
leads to a contradiction when $f_1\in \Omega$ since $\ker \chi\Pi_{\la_0}|_{C^\infty(K)}\cap \Omega=\emptyset$ and so
$\chi R(\la)f$ has a singularity at $\la=\la_0$. 
We thus obtain our conclusion. The same method applies when $X$ is Schottky and the finer estimates are valid.
\qed

\section{Conformal resonances}
In this section, we try to explain the special cases $\delta\in n/2-\nn$ in term of conformal theory
of the conformal infinity.
As emphasized before, a convex co-compact hyperbolic manifold $(X,g)$
compactifies into a smooth compact manifold with boundary $\bar{X}=X\cup \pl\bar{X}$,
where $\pl\bar{X}=\Gamma\backslash\Omega$ if $\Omega$ is the domain of discontinuity 
of the group $\Gamma$ defined in the introduction.
If $x$ is a smooth boundary defining function of $\pl\bar{X}$, $x^2g$ extends smoothly to $\bar{X}$ as a metric,
the restriction 
\[h_0=x^2g|_{T\pl\bar{X}}\]
is a metric on $\pl\bar{X}$ inherited from $g$ but depending on the choice of $x$, however 
its conformal class $[h_0]$ is clearly independent of $x$, it is then called the \emph{conformal infinity} of $X$. 
By Graham-Lee \cite{GRL,GR}, there is an identification between a particular class of boundary defining functions
and elements of the class $[h_0]$: indeed, for any $h_0\in[h_0]$, there exists near $\pl\bar{X}$ a unique boundary
defining function $x$ such that $|dx|_{x^2g}=1$ and $x^2g|_{T\pl\bar{X}}=h_0$, this function will be called a \emph{geodesic boundary defining function}.\\

We now recall the definition of the scattering operator $S(\la)$ as in \cite{GRZ,JSB}. Let $\la\in\cc$ with $\Re(\la)\notin n/2+\zz$ and let $x$ be a geodesic boundary defining function, then for all $f\in C^{\infty}(\pl\bar{X})$
there exists a unique function $F(\la,f)\in C^\infty(X)$ which satisfies 
the boundary value problem
\[
\left\{\begin{array}{l}
(\Delta_X-\la(n-\la))F(\la,f)=0,\\ 
\exists F_1(\la,f),F_2(\la,f)\in C^\infty(\bar{X}) \textrm{ such that }\\
F(\la,f)=x^{n-\la}F_1(\la,f)+x^\la F_2(\la,f) \textrm{ and }  F_1(\la,f)|_{\pl\bar{X}}=f.
\end{array}\right.\] 
Then the operator $S(\la):C^\infty(\pl\bar{X})\to C^{\infty}(\pl\bar{X})$ is defined by 
\[S(\la)f =F_2(\la,f)|_{\pl\bar{X}}.\]
It is clear that $S(\la)$ depends on choice of $x$, but it is conformally covariant under change of boundary 
defining function: if $\hat{x}:=xe^{\omega}$ is another such function, then the related scattering operator is
\[\hat{S}(\la)=e^{-\la\omega_0}S(\la)e^{(n-\la)\omega_0}, \quad \omega_0:=\omega|_{\pl\bar{X}}.\]
It is proved in \cite{GRZ} that $S(\la)$ has simple poles at $\la=n/2+k$ for all $k\in\nn$, and  after renormalizing $S(\la)$ into 
\[\mc{S}(\la):=2^{2\la-n}\frac{\Gamma(\la-\ndemi)}{\Gamma(\ndemi-\la)}S(\la)\]
we obtain by the main result of \cite{GRZ} that $\mc{S}(n/2+k)=P_k$ is the k-th GJMS conformal Laplacian
on $(\pl\bar{X},h_0)$ defined previously in \cite{GJMS}. 
In general $\mc{S}(\la)$ is a pseudodifferential operator of order $2\la-n$ 
with principal symbol $|\xi|_{h_0}^{2\la-n}$ but for $\la=n/2+k$, it becomes 
differential.

\begin{prop}\label{confresonance}
If $\delta=n/2-k$ with $k\in\nn$, then the $j$-th GJMS conformal Laplacian
$P_j>0$ for $j<k$ while $P_k$ has a kernel of dimension $1$ with 
eigenvector given by $f_{n/2-k}$ defined below in \eqref{fndemik} in term of Patterson-Sullivan measure. 
\end{prop}
\textsl{Proof}: Let us fix $\delta\in(0,n/2)$ not necessarily in $n/2-\nn$ for the moment.
In \cite{G}, the first author studied the relation between poles of resolvent and poles of scattering operator.
If $\la\in\cc$, we define its resonance multiplicity by 
\[m(\la):=\rang \Big(\textrm{Res}_{s=\la}((2s-n)R(s))\Big)\]
while its scattering pole multiplicity is defined by 
\[\nu(\la):=-\tra \Big(\textrm{Res}_{s=\la}(\pl_s\mc{S}(s)\mc{S}^{-1}(s))\Big).\]
We proved in \cite{G} (see also \cite{GN} for point in pure point spectrum) that for $\Re(\la)<n/2$
\[\nu(\la)=m(\la)-m(n-\la)+\indic_{\ndemi-\nn}(\la)\dim\ker \mc{S}(n-\la),\]   
which in our case reduces to 
\begin{equation}\label{reduce}
\nu(\la)=m(\la)+\indic_{\ndemi-\nn}(\la)\dim\ker \mc{S}(n-\la)
\end{equation}
by the holomorphy of $R(\la)$ in $\{\Re(\la)\geq n/2\}$, stated in Proposition \ref{pa}. 
We know from \cite{JSB,GRZ} that the Schwartz kernel of $\mc{S}(\la)$ is related to that of $R(\la)$ by 
\begin{equation}\label{noyaus}
\mc{S}(\la;y,y')=2^{2\la-n+1}\frac{\Gamma(\la-\ndemi+1)}{\Gamma(\ndemi-\la)}[x^{-\la}{x'}^{-\la}R(\la;x,y,x'y')]|_{x=x'=0}
\end{equation}
where $(x,y)\in [0,\eps)\x\pl\bar{X}$ are coordinates in a collar neighbourhood of $\pl\bar{X}$,
$x$ being the geodesic boundary defining function used to define $\mc{S}(\la)$. This implies with 
Proposition \ref{pa} that $\mc{S}(\la)$ is analytic in $\{\Re(\la)>\delta\}$ and 
has a simple pole at $\delta$ with residue
\[\textrm{Res}_{\la=\delta}\mc{S}(\la)=A_X \frac{2^{-2k+1}}{(k-1)!}f_\delta\otimes f_\delta,
\quad f_\delta:=(x^{-\delta}u_\delta)|_{x=0}.\]
Note that Perry \cite{P2} proved that $f_\delta$ is well defined and in $C^\infty(\pl\bar{X})$.
The functional equation 
$\mc{S}(\la)\mc{S}(n-\la)=\rm{Id}$ (see for instance Section 3 of \cite{GRZ}) 
and the fact that $\mc{S}(\la)$ is analytic in $\{\Re(\la)>\delta\}$ clearly imply 
that $\ker\mc{S}(\la)=0$ for $\Re(\la)\in(\delta,n-\delta)$, thus in particular 
$\ker P_j=0$ for any $j\in \nn$ with $j<n/2-\delta$. Moreover, using \cite[Lemma 4.16]{PP} and the fact that
$m_{n/2}=0$ since $R(\la)$ is holomorphic in $\{\Re(\la)>\delta\}$, one obtains $S(n/2)=\textrm{Id}$ thus
$\mc{S}(\la)>0$ for all $\la\in (\delta,n-\delta)$ by continuity of $\mc{S}(\la)$ with respect to $\la$.
We also deduce from the functional equation and the holomorphy of $\mc{S}(s)$ at $n-\delta$ that 
\[\mc{S}(n-\delta)f_{\delta}=0.\]
We thus see from this discussion and Proposition \ref{pa} that, in \eqref{reduce}, the relation 
$m(\delta)=\nu(\delta)=1$ holds when $\delta\notin n/2-\nn$ while 
$\nu(\delta)=\dim\ker P_k$ when $\delta=n/2-k$ with $k\in\nn$ since $m(\delta)=0$ in that case by holomorphy
of $R(\la)$ at $\delta=n/2-k$. To compute $\dim\ker P_k$ when $\delta=n/2-k$, one can use for instance
Selberg's zeta function. Indeed by Proposition 2.1 of \cite{P2}, $Z(\la)$ has a simple zero
at $\delta$ but  it follows from Theorems 1.5-1.6 of Patterson-Perry \cite{PP} 
that $Z(\la)$ has a zero at $\la=n/2-k$ of order $\nu(n/2-k)$ if $k\in\nn, k<n/2$, therefore $\nu(n/2-k)=1$ and thus 
\[\dim\ker P_k=1.\]
One can now describe a bit more precisely the function $f_\delta$.
The Poisson kernel of Proposition \ref{pa} in the half-space model $\rr_y^{n}\x\rr_{y_{n+1}}^+$ of 
$\hh^{n+1}$ is 
\[\mc{P}(\la; y,y_{n+1},y')=\frac{y_{n+1}}{y_{n+1}^2+|y-y'|^2}\]
thus if $x$ is the boundary defining function used to define $\mc{S}(\la)$ and if $(\pi_\Gamma^*x/y_{n+1})
|_{y_{n+1}=0}=k(y)$ (recall $\pi_\Gamma,\bar{\pi}_\Gamma$ are the projections of (\ref{pigamma})) 
for some $k(y)\in C^\infty(\rr^n)$, so we  can describe rather explicitely $f_\delta$, we have
\begin{equation}\label{fndemik}
\bar{\pi}_\Gamma^*f_{\delta}(y)=k(y)^{-\delta}\int_{\rr^n}|y-y'|^{-2\delta}d\mu_\Gamma(y'), \quad y\in \Omega.
\end{equation}
\qed\\
 
To summarize the discussion, if $\delta<n/2$, the Patterson function $u_\delta$ is an eigenfunction for $\Delta_X$ with eigenvalue $\delta(n-\delta)$, it is not an $L^2$ eigenfunction though and it has leading asymptotic behaviour $u_\delta\sim x^{\delta}f_\delta$ as $x\to 0$, where $f_\delta\in C^\infty(\pl\bar{X})$
is in the kernel of the boundary operator $\mc{S}(n-\la)$. When $\delta\notin n/2-\nn$, this is a resonant state 
for $\Delta_X$ with associated resonance $\delta$ while when $\delta\in n/2-\nn$ it is still a generalized eigenfunction  
of $\Delta_X$ but not a resonant state anymore, and $\delta$ is not a resonance yet in that case: the resonance  
disappear when $\delta$ reaches $n/2-k$ and instead the $k$-th GJMS at $\pl\bar{X}$ gains an element in its kernel given by the leading coefficient of $u_{n/2-k}$ in the asymptotic at the boundary.\\ 

\textsl{Remark}: Notice that the positivity of $P_j$ for $j<n/2-\delta$ has been proved by Qing-Raske \cite{QR}
and assuming a positivity of Yamabe invariant of the boundary. Our proof allows to remove the assumption on the Yamabe 
invariant, which, as we showed, is automatically satisfied if $\delta<n/2$.

\end{document}